# Approximate method for cutting pattern optimization of frame-supported and pneumatic membrane structures

Makoto Ohsaki[1*], Jun Fujiwara[2] and Fumiyoshi Takeda[3]

[1]Department of Architecture and Architectural Engineering, Kyoto University,
Kyoto-Daigaku Katsura, Nishikyo, Kyoto 615-8540, Japan, ohsaki@archi.kyoto-u.ac.jp

[2]National Research Institute for Earth Science and Disaster Resilience, Japan, j.fujiwara@bosai.go.jp

[3]Taiyo Kogyo Corporation, Japan, tf001163@mb.taiyokogyo.co.jp

**Abstract**

A computationally efficient method is presented for approximate optimization of cutting pattern of frame-supported and pneumatic membrane structures. The plane cutting sheet is generated by minimizing the error from the shape obtained by reducing the stress from the desired curved shape. The equilibrium shape is obtained solving a minimization problem of total strain energy. The external work done by the pressure is also incorporated for analysis of pneumatic membrane. An approximate method is also proposed for analysis of an Ethylene TetraFluoroEthylene (ETFE) film, where elasto-plastic behavior is modeled as a nonlinear elastic material under monotonic loading condition. Efficiency of the proposed method is demonstrated through examples of a frame-supported PolyVinyl Chloride (PVC) membrane structure and an air pressured square ETFE film.

## 1. Introduction

In the field of civil and architectural engineering, curved surface structures, which are formed with pre-tensioned coated fabrics or synthetic resin films, are widely used for lightweight and long-span roofs or facades of stadiums, gymnasium and so on [1]. This kind of structure called membrane structure is also used for solar panels in space [2]. In the first step of designing process of a membrane structure, a target surface is found under given stress distribution and boundary conditions, where in most of cases the target stress distribution is assumed to be uniform tension. Since coated fabrics and resin films are fabricated as planar sheets, the target surface is realized by connecting planar cutting sheets and clamping their edges to boundary structures such as steel frames and cables. Since curved surfaces are not generally developable to a plane, there must be some errors between the target stress distribution and the stresses of installed membrane structure.

In the process of designing a cutting pattern, it is important to achieve a moderately uniform stress distribution to prevent fracture and slackening. Although there are many methods for cutting pattern optimization, most of the methods should carry out finite element analysis many times [3]. Ohsaki and Fujiwara [4], Ohsaki and Uetani [5] and Uetani *et al*. [6] proposed a method for inversely generating the shape of cutting pattern by removing the stresses from the self-equilibrium shape. Kim and Lee [7] presented an approximate method for a membrane divided into many strips. Punurai *et al*. [8] used genetic algorithm for cutting pattern optimization. Philipp *et al*. [9] proposed an updated reference strategy for form-finding of membrane structures supported by bending-active members. Extensive researches have also been carried out for surface flattening in the field of computer aided design [10, 11].

The woven fabric membrane material has orthotropic property [12, 13], which is often simplified to isotropic or orthotropic linear elastic material. Therefore, the equilibrium shape analysis for specified cutting pattern can be formulated as a forced displacement analysis problem, which can be solved by geometrically nonlinear analysis [14, 15], dynamic relaxation method [16], or minimization of the total strain energy [17, 18]. However, there exists difficulty for form-finding and equilibrium shape analysis of pneumatic membrane structures that are formed and strengthened by applying air pressure [14, 19]. Bouzidi and Le Van [18] proposed an energy minimization method for analysis of pneumatic membrane structures supported by rigid boundary incorporating a pressure potential function. However, Ethylene TetraFluoroEthylene (ETFE) film, which is often used for pneumatic membrane structures, has an elasto-plastic property; therefore, it is difficult to optimize the shape of the surface using a gradient-based optimization algorithm. Dinh *et al*. [20] proposed an elasto-plastic material model for equilibrium shape analysis, which is applied to shape optimization of membrane with boundary cables [21].

In this paper, we present a computationally efficient iterative method for approximate optimization of cutting pattern of membrane structures. The plane cutting sheet is generated by minimizing the error from the shape obtained by reducing the stress from the desired curved shape, which is discretized into triangular finite elements. The equilibrium shape corresponding to the specified cutting pattern is obtained by energy minimization. The external work done by the air pressure is also incorporated for analysis of pneumatic membrane structures. Efficiency of the proposed method is demonstrated through examples of a frame-supported PolyVinyl Chloride (PVC) membrane structure and an air-pressured square ETFE film.

## 2. Energy minimization for equilibrium shape analysis

Self-equilibrium shape of a frame-supported membrane structure is generated by connecting plane cutting sheets to the boundary. This process is regarded as a forced deformation analysis problem, which can be solved by minimizing an energy function. In this section, we briefly present the formulation of strain energy minimization, and extend it to pneumatic membrane structures [22].

**2.1 Strain energy minimization for frame-supported membrane**

Consider a curved membrane structure discretized into triangular finite elements with constant stress in plane stress state. The material property is assumed to be linear orthotropic elastic. Let $E_x$ and $E_y$ denote the elastic moduli in the two principal directions $(x_p, y_p)$ of the material. The shear modulus is denoted by $G$. The following notations are used:

$$\beta = \frac{E_x}{E_y}, \quad \kappa = \frac{G}{E_y} \tag{1}$$

Then the constitutive matrix $\mathbf{D} \in \mathbb{R}^{3\times 3}$ defining the orthotropic elastic material property of the membrane is given as

$$\mathbf{D} = \frac{E_y}{1-\beta v_{xy}^2} \begin{bmatrix} \beta & \beta v_{xy} & 0 \\ \beta v_{xy} & 1 & 0 \\ 0 & 0 & \kappa(1-\gamma v_{xy}^2) \end{bmatrix} \tag{2}$$

In the following, stress is evaluated as the force per unit length of the section. The vectors of inplane stress and strain with respect to the principal coordinates of element $k$ are denoted by $\boldsymbol{\sigma}_k = (\sigma_{kx}, \sigma_{ky}, \tau_k)^T$ and $\boldsymbol{\varepsilon}_k = (\varepsilon_{kx}, \varepsilon_{ky}, \gamma_k)^T$, respectively, where $\sigma_{kx}$ and $\sigma_{ky}$ are the stresses in the directions of $x_p$ and $y_p$; $\varepsilon_{kx}$ and $\varepsilon_{ky}$ are the strains in the directions of $x_p$ and $y_p$; $\tau_k$ is the shear stress; and $\gamma_k$ is the shear strain. The stress-strain relation of the $k$th element is written as

$$\boldsymbol{\sigma}_k = \mathbf{D}\boldsymbol{\varepsilon}_k \tag{3}$$

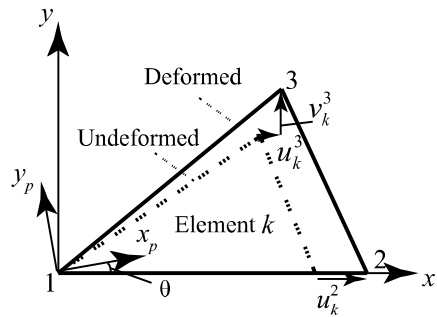

Figure 1: Local coordinates, node numbers, principal directions, and deformation

of triangular element $k$.

Figure 1 defines the local (*x*,*y*)-coordinates and relative displacements from undeformed state to deformed state of element *k* consisting of nodes 1, 2, and 3. The displacements are fixed in *x*- and *y*-directions at node 1, and fixed in *y*-direction at node 2. Let $u_k^i$ and $v_k^i$ denote the relative displacements of node *i* from node 1 of element *k*, which are combined to a vector $\mathbf{u}_k = (u_k^2, u_k^3, v_k^3)^T$ after constraining the rigid body displacements and rotation. The strain-displacement relation is written using a strain-displacement relation matrix $\mathbf{C} \in \mathbb{R}^{3 \times 3}$ and the rotation matrix $\mathbf{R} \in \mathbb{R}^{3 \times 3}$ from the local to the principal coordinate systems as

$$\boldsymbol{\varepsilon}_k = \mathbf{RCu}_k \qquad (4)$$

See a textbook of finite element analysis for details.

The vector consisting of the global $(X, Y, Z)$-coordinates of the nodes on the surface is denoted by $\mathbf{X} \in \mathbb{R}^{3n}$, where *n* is the number of nodes. In the process of finding the equilibrium shape using forced displacement analysis, the strain energy $S(\mathbf{X})$ of the membrane is regarded as a function of $\mathbf{X}$ as

$$S(\mathbf{X}) = \frac{1}{2} \sum_{k=1}^{m} A_k \boldsymbol{\varepsilon}_k(\mathbf{X})^T \mathbf{D} \boldsymbol{\varepsilon}_k(\mathbf{X}) \qquad (5)$$

where *m* is the number of elements, and $A_k(\mathbf{X})$ is the area of the *k*th element.

Assuming that the self-weight is negligibly small compared to the tension forces in the membrane, the equilibrium shape is obtained by solving the following problem of minimizing $S(\mathbf{X})$ under an appropriate boundary conditions restricting rigid body motions:

$$\begin{aligned} \text{Minimize} \quad & S(\mathbf{X}) \\ \text{subject to} \quad & \mathbf{X}_\text{L} \leq \mathbf{X} \leq \mathbf{X}_\text{U} \end{aligned} \qquad (6)$$

where $\mathbf{X}_\text{L}$ and $\mathbf{X}_\text{U}$ are the lower and upper bounds for $\mathbf{X}$.

**2.1 Energy minimization for pneumatic membrane**

For a frame-supported pneumatic membrane structure with air pressure *p*, the pressure potential energy $W(\mathbf{X})$ is given as [14, 23]

$$W(\mathbf{X}) = pV(\mathbf{X}) \qquad (7)$$

where $V(\mathbf{X})$ is the volume of the air contained by the membrane structure. Let $\mathbf{X}_i^\text{S} = (X_{i1}^\text{S}, X_{i2}^\text{S}, X_{i3}^\text{S})^T$ denote the coordinate vector of node *i* on the curved surface, where the subscripts 1, 2, and 3 correspond to *X*-, *Y*-, and *Z*-components; hence, $\mathbf{X}$ consists of the free components of $\mathbf{X}_i^\text{S}$. The coordinate vector $\overline{\mathbf{X}}_k^\text{S} \in \mathbb{R}^3$ of the center of gravity of element *k* is defined as

$$\overline{\mathbf{X}}_k^{\mathrm{S}}(\mathbf{X}) = \frac{1}{3}\sum_{i \in I_k} \mathbf{X}_i^{\mathrm{S}} \tag{8}$$

where $I_k$ is the set of indices of three nodes of element $k$.

The functions $W(\mathbf{X})$ and $V(\mathbf{X})$ are expressed as summation of element values $W_k^{\mathrm{e}}(\mathbf{X})$ and $V_k^{\mathrm{e}}(\mathbf{X})$, respectively. Variation of $W(\mathbf{X})$ in the normal direction of surface is computed as

$$\begin{aligned}
\delta W(\mathbf{X}) &= \sum_{k=1}^{m} \delta W_k^{\mathrm{e}}(\mathbf{X}) \\
&= p \sum_{k=1}^{m} \delta V_k^{\mathrm{e}}(\mathbf{X}) \\
&= p \sum_{k=1}^{m} \sum_{i \in I_k} \sum_{j=1}^{3} \left( \frac{\partial V_k^{\mathrm{e}}(\mathbf{X})}{\partial X_{ij}^{\mathrm{S}}} \delta X_{ij}^{\mathrm{S}} \right)
\end{aligned} \tag{9}$$

where $X_{ij}^{\mathrm{S}}$ is the $j$th component of $\mathbf{X}_i^{\mathrm{S}}$. Let $\mathbf{n}_k(\mathbf{X}) \in \mathbb{R}^3$ denote the unit normal vector of the $k$th triangular element, and suppose node $i$ belongs to element $k$. Then the differential coefficient of $V_k^{\mathrm{e}}(X_{ij}^{\mathrm{S}})$ with respect to $X_{ij}^{\mathrm{S}}$ in the normal direction of element $k$ is written as

$$\frac{\partial V_k^{\mathrm{e}}(\mathbf{X})}{\partial X_{ij}^{\mathrm{S}}} = A_k \mathbf{n}_k(\mathbf{X})^T \frac{\partial \overline{\mathbf{X}}_k^{\mathrm{S}}(\mathbf{X})}{\partial X_{ij}^{\mathrm{S}}} \tag{10}$$

Using Eq. (10), $\delta W(\mathbf{X})$ is further modified as

$$\begin{aligned}
\delta W(\mathbf{X}) &= p \sum_{k=1}^{m} \sum_{i \in I_k} \sum_{j=1}^{3} \left[ A_k \mathbf{n}_k(\mathbf{X})^T \left( \frac{\partial \overline{\mathbf{X}}_k^{\mathrm{S}}(\mathbf{X})}{\partial X_{ij}^{\mathrm{S}}} \delta X_{ij}^{\mathrm{S}} \right) \right] \\
&= p \sum_{k=1}^{m} \sum_{i \in I_k} \sum_{j=1}^{3} \left[ \frac{A_k}{3} \mathbf{n}_k(\mathbf{X})^T \left( \frac{\partial \mathbf{X}_i^{\mathrm{S}}(\mathbf{X})}{\partial X_{ij}^{\mathrm{S}}} \delta X_{ij}^{\mathrm{S}} \right) \right] \\
&= p \sum_{k=1}^{m} \sum_{i \in I_k} \sum_{j=1}^{3} \left\{ \frac{A_k}{3} \mathbf{n}_k(\mathbf{X})^T \left[ \begin{pmatrix} 1 \\ 0 \\ 0 \end{pmatrix} \delta X_{i1}^{\mathrm{S}} + \begin{pmatrix} 0 \\ 1 \\ 0 \end{pmatrix} \delta X_{i2}^{\mathrm{S}} + \begin{pmatrix} 0 \\ 0 \\ 1 \end{pmatrix} \delta X_{i3}^{\mathrm{S}} \right] \right\} \\
&= p \sum_{k=1}^{m} \sum_{i \in I_k} \left( \frac{A_k}{3} \mathbf{n}_k(\mathbf{X})^T \delta \mathbf{X}_i^{\mathrm{S}} \right) \\
&= p \sum_{i=1}^{n} \sum_{k \in K_i} \left( \frac{A_k}{3} \mathbf{n}_k(\mathbf{X})^T \delta \mathbf{X}_i^{\mathrm{S}} \right)
\end{aligned} \tag{11}$$

where $K_i$ is the set of elements connected to node $i$.

Using the divergence theorem, Bouzidi and Le Van [14] derived the following expression $W^*(\mathbf{X})$ for a membrane discretized by triangular finite elements:

$$W^*(\mathbf{X}) = p \int_V dV$$
$$= \frac{p}{3} \int_V \mathrm{div}\, \mathbf{X} dV$$
$$= \frac{p}{3} \int_\Gamma \mathbf{X}^T \mathbf{n}(\mathbf{X}) d\Gamma \quad (12)$$
$$= \frac{p}{3} \sum_{k=1}^m A_k \mathbf{n}_k(\mathbf{X})^T \overline{\mathbf{X}}_k^S(\mathbf{X})$$

where $\Gamma$ is the boundary of the space covered by the membrane. However, the term 1/3 in the right-hand-side of Eq. (12) is not necessary as confirmed below, because we should consider shape variation only in the normal direction of surface. Actually they did not use 1/3 in the numerical examples. Therefore, we use the following definition:

$$W(\mathbf{X}) = p \sum_{k=1}^m A_k \mathbf{n}_k(\mathbf{X})^T \overline{\mathbf{X}}_k^S(\mathbf{X}) \quad (13)$$

The equilibrium shape is obtained by minimizing the total potential energy $\Pi(\mathbf{X})$ defined as

$$\Pi(\mathbf{X}) = S(\mathbf{X}) - W(\mathbf{X}) \quad (14)$$

Differentiation of $\Pi(\mathbf{X})$ with respect to $X_{ij}^S$ leads to

$$\frac{\partial \Pi(\mathbf{X})}{\partial X_{ij}^S} = \sum_{k \in K_i} \left( A_k \boldsymbol{\varepsilon}_k(\mathbf{X})^T \mathbf{D} \frac{\partial \boldsymbol{\varepsilon}_k(\mathbf{X})}{\partial X_{ij}^S} \right) - p \sum_{k \in K_i} \left( A_k \mathbf{n}_k(\mathbf{X})^T \frac{\partial \overline{\mathbf{X}}_k^S(\mathbf{X})}{\partial X_{ij}^S} \right)$$
$$- p \sum_{k \in K_i} \left( A_k \frac{\partial \mathbf{n}_k(\mathbf{X})^T}{\partial X_{ij}^S} \overline{\mathbf{X}}_k^S(\mathbf{X}) \right) \quad (15)$$
$$= \sum_{k \in K_i} \left( A_k \boldsymbol{\varepsilon}_k^T(\mathbf{X}) \mathbf{D} \frac{\partial \boldsymbol{\varepsilon}_k(\mathbf{X})}{\partial X_{ij}^S} \right) - \frac{p}{3} \sum_{k \in K_i} \left( A_k n_{kj}(\mathbf{X}) \right) - p \sum_{k \in K_i} \left( A_k \frac{\partial \mathbf{n}_k(\mathbf{X})^T}{\partial X_{ij}^e} \overline{\mathbf{X}}_k^S(\mathbf{X}) \right)$$

where Eq. (8) has been used, and $n_{kj}$ is the $j$th component of $\mathbf{n}_k$. By contrast, the equilibrium equation is given as

$$\sum_{k \in K_i} \left( A_k \boldsymbol{\varepsilon}_k^T(\mathbf{X}) \mathbf{D} \frac{\partial \boldsymbol{\varepsilon}_k(\mathbf{X})}{\partial X_{ij}^S} \right) - \frac{p}{3} \sum_{k \in K_i} \left( A_k n_{kj}(\mathbf{X}) \right) = 0, \quad (i=1,\ldots,n; j=1,2,3) \quad (16)$$

Therefore, the third term, denoted by $e_{ij}$, in the right-hand-side of Eq. (15) remains as an additional term, which is rewritten as

$$e_{ij} = -p \sum_{k \in K_i} \left( A_k \frac{\partial \mathbf{n}_k(\mathbf{X})^T}{\partial X_{ij}^{\mathrm{e}}} \bar{\mathbf{X}}_k^{\mathrm{S}}(\mathbf{X}) \right) \quad (17)$$

It is easily seen that the direction of the normal vector does not change due to the variation of nodal coordinates in the direction of the plane of the element. Furthermore, the variation of normal vector due to the variation of nodes in the normal direction can be neglected within the first-order approximation. Therefore, the term $\partial \mathbf{n}_k(\mathbf{X})^T / \partial X_{ij}^{\mathrm{S}}$ vanishes, and $e_{ij} = 0$ is satisfied; accordingly, the equilibrium shape of pneumatic membrane can be found by solving the following problem of minimizing the energy function $\Pi(\mathbf{X})$:

$$\begin{aligned} \text{Minimize} &\quad \Pi(\mathbf{X}) \\ \text{subject to} &\quad \mathbf{X}_{\mathrm{L}} \leq \mathbf{X} \leq \mathbf{X}_{\mathrm{U}} \end{aligned} \quad (18)$$

## 3. Approximate optimization of cutting pattern

Suppose a target shape of membrane surface is given by a designer. Since the arbitrary defined shape cannot be realized by connecting plane cutting sheets, a simple update rule of the stress parameters called *reduction stress* is proposed for approximate cutting pattern optimization. The objective here is to find the shape of cutting pattern for minimizing the variation of stresses from the *ideal target stresses*, which usually represent uniform tension state. The method is based on the inverse process of generating a plane sheet from a curved surface by reducing the stress [1].

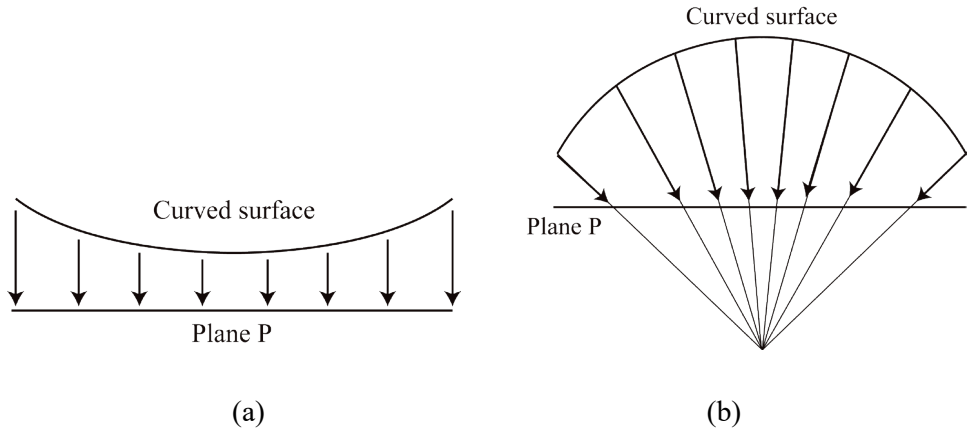

Figure 2: Projection of surface to a plane; (a) small curvature, (b) large curvature.

We assign a plane P near the target surface, and project the triangular mesh onto P to generate the triangular mesh of the cutting sheet as an initial solution for optimization. If the curvature of surface is small, the plane P may be parallel to a tangent plane at the center of surface as illustrated in Fig. 2(a). By contrast, if the curvature is large as illustrated in Fig. 2(b), the surface may be projected in the direction of a point near the center of curvature.

The vector consisting of $(X^P, Y^P)$-coordinates of the nodes of the cutting sheet on plane P is denoted by $\mathbf{X}^P \in \mathbb{R}^{2n}$. Let $L_{k1}^P(\mathbf{X}^P)$, $L_{k2}^P(\mathbf{X}^P)$, and $L_{k3}^P(\mathbf{X}^P)$ denote the edge lengths of the $m$ triangular elements on P, which are functions of $\mathbf{X}^P$. Let $L_{k1}^0$, $L_{k2}^0$, and $L_{k3}^0$ denote the unstressed lengths of three edges of element $k$ after removing the stress from the triangular element on the curved surface. Then the following problem is solved to minimize the sum of error between the lengths of a pair of edges on the surface and on the plane:

$$\text{Minimize} \quad F(\mathbf{X}^P) = \sum_{k=1}^{m} \sum_{i=1}^{3} \kappa_{ki} (L_{ki}^P(\mathbf{X}^P) - L_{ki}^0)^2 \quad (19)$$
$$\text{subject to} \quad \mathbf{X}_L^P \leq \mathbf{X}^P \leq \mathbf{X}_U^P$$

where $\kappa_{ki}$ is a weight parameter, and $\mathbf{X}_L^P$ and $\mathbf{X}_U^P$ are the lower and upper bounds for $\mathbf{X}^P$. In the following examples a large weight is given for a short edge as $\kappa_{ki} = 1/L_{ki}^0$ to prevent reversal of a triangle on the plane.

Outline of the algorithm is illustrated in Fig. 3 and summarized as follows:

Step 1: Assign the target shape of curved surface, boundary conditions, and generate triangular meshes on the target surface. Compute the edge lengths $L_{k1}$, $L_{k2}$, and $L_{k3}$ of the $k$th triangle of the target surface. Specify the ideal target stresses $\sigma_{k1}^*$, $\sigma_{k2}^*$, and $\tau_k^* (=0)$, and initialize the step counter $s = 0$.

Step 2: Project the triangulated surface to the plane P to find the initial values of $\mathbf{X}^P$ for solving the optimization problem (19).

Step 3: Specify the reduction stresses $\hat{\sigma}_{k1}^s$, $\hat{\sigma}_{k2}^s$, and $\hat{\tau}_k^s (=0)$ in principal directions of material. Remove the stress from the triangular elements on the curved surface, and compute the unstressed edge lengths $L_{k1}^0$, $L_{k2}^0$, and $L_{k3}^0$.

Step 4: Find the nodal coordinates of the cutting sheet on plane P by solving problem (19).

Step 5: Carry out equilibrium shape analysis by solving the energy minimization problems (6) and (18), respectively, for frame-supported membrane and pneumatic membrane to find the nodal coordinates on surface and the stresses $\sigma_{k1}^s$, $\sigma_{k2}^s$, and $\tau_k^s$ at the self-equilibrium state.

Step 6: Let $s \leftarrow s+1$, and modify the reduction stresses $\hat{\sigma}_{k1}^s$ and $\hat{\sigma}_{k2}^s$ as

$$\hat{\sigma}_{k1}^{s+1} = \hat{\sigma}_{k1}^s + c(\sigma_{k1}^* - \sigma_{k1}^s), \quad \hat{\sigma}_{k2}^{s+1} = \hat{\sigma}_{k2}^s + c(\sigma_{k2}^* - \sigma_{k2}^s), \quad \hat{\tau}_k^{s+1} (=0) \qquad (20)$$

where $c$ is the parameter for adjusting the convergence property. Also replace the target surface with the equilibrium surface obtained in Step 5.

Step 7: Go to Step 2, if termination condition is not satisfied.

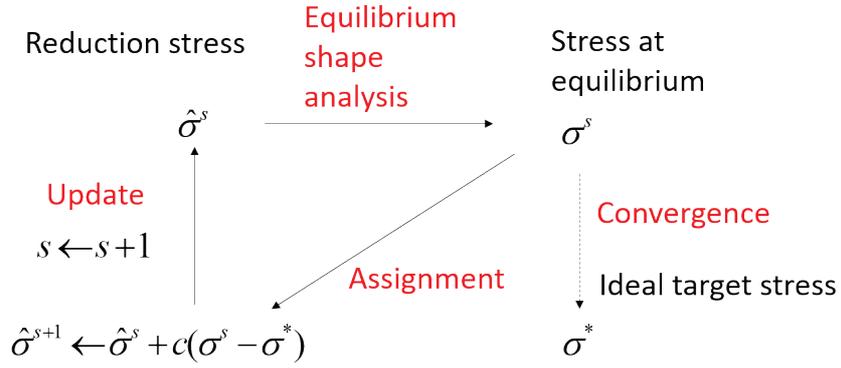

Figure 3: Scheme of approximate cutting pattern optimization.

This process is illustrated using a single cable that is pin-supported at the both ends. The units are omitted for brevity, and Young's modulus is assumed to be 1. The span length at equilibrium is $L=1.2$. The problem is to find the unstressed length of the cable so that the stress after connecting the both ends to the supports becomes 0.1. In order to follow the process for membrane structures, for which the shape of a triangular element at equilibrium is not known, the initial guess for the length at equilibrium is assumed to be $\bar{L}=1.0$, although it is known as 1.2. Parameter $c$ is equal to 1.0. The process is listed as follows:

Step 1. Assign the reduction stress $\hat{\sigma}^1$ equal to the ideal target stress as $\hat{\sigma}^1 = \sigma^* = 0.1$.
Compute the unstressed length as $L_0 = \bar{L} \times (1.0 - 0.1) = 1.0 \times (1.0 - 0.1) = 0.9$.
Compute the stress at equilibrium as $\sigma^1 = (L - L_0)/L_0 = (1.2 - 0.9)/0.9 = 0.33$.

Step 2. Update the reduction stress as
$\hat{\sigma}^2 = \hat{\sigma}^1 + c(\sigma^* - \sigma^1) = 0.1 + 1.0 \times (0.1 - 0.33) = -0.13$.
Compute the unstressed length as $L_0 = \bar{L} \times (1.0 + 0.13) = 1.0 \times (1.0 + 0.13) = 1.13$.
Compute the stress at equilibrium as $\sigma^2 = (L - L_0)/L_0 = (1.2 - 1.13)/1.13 = 0.06$.

Step 3. Update the reduction stress as
$\hat{\sigma}^3 = \hat{\sigma}^2 + c(\sigma^* - \sigma^2) = -0.13 + 1.0 \times (0.1 - 0.06) = -0.09$.
Compute the unstressed length as $L_0 = \bar{L} \times (1.0 + 0.09) = 1.0 \times (1.0 + 0.09) = 1.09$.

Compute the stress at equilibrium as $\sigma^2 = (L - L_0)/L_0 = (1.2 - 1.09)/1.09 = 0.10$.

This way, the correct unstressed length has been found with three iterations.

As seen from this result of a cable, initial guess of the length $\bar{L}$ does not have any effect on the final shape. Accordingly, the reduction stress is a controlling parameter, and is not related to the final stress at equilibrium. Therefore, for a membrane structure, the accuracy of the projected shape on the plane is not very important.

## 4. Analysis of ETFE film

ETFE film is usually modeled as elasto-plastic material with von Mises yield criterion [24, 25]. However, the method of cutting pattern design described in Sec. 3 is applicable only to elastic membrane material. Furthermore, only monotonic loading is considered in the process of increasing the pressure to reach the equilibrium shape. Therefore, in this section, an approximation method is presented for ETFE film using nonlinear elastic material model.

The relation between stress and strain in uniaxial tension is often modeled as bilinear relation, which is identified by experiments as shown in Fig. 4, where $\sigma^Y$ and $\varepsilon^Y$ are the yield stress and strain. For a continuum with von Mises criterion, Fig. 4 may be regarded as a relation between equivalent stress and strain. Since the stiffness after yielding is smaller than the initial elastic stiffness, almost uniform stress distribution can be expected, if the target stress $\sigma^*$ is larger than $\sigma^Y$.

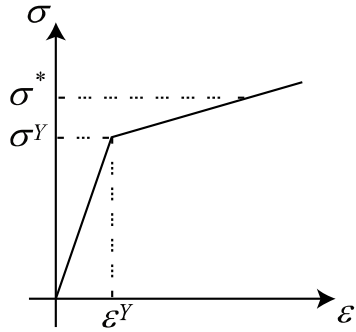

Figure 4: Relation between stress $\sigma$ and strain $\varepsilon$ of ETFE sheet under uniaxial tension.

Let $E$, $H$ and $\mathbf{D}_1$ denote Young's modulus, hardening coefficient and elastic constitutive matrix, respectively. To model the stress-strain relation as nonlinear elastic, the constitutive matrix after yielding is given as

$$\mathbf{D}_2 = \frac{H}{E}\mathbf{D}_1 \tag{21}$$

In the process of energy minimization for finding the self-equilibrium shape, the trial stress vector $\tilde{\boldsymbol{\sigma}} = (\tilde{\sigma}_x, \tilde{\sigma}_y, \tilde{\tau})^T$ with respect to the element coordinates $(x, y)$ is first calculated as follows from the strain vector $\boldsymbol{\varepsilon} = (\varepsilon_x, \varepsilon_y, \tau)^T$ with respect to the $(x, y)$-coordinates, which is calculated from the nodal coordinates on the surface:

$$\tilde{\boldsymbol{\sigma}} = \mathbf{D}_1 \boldsymbol{\varepsilon} \tag{22}$$

Then, the equivalent stress $\tilde{\sigma}_{eq}$ corresponding to $\tilde{\boldsymbol{\sigma}}$ is obtained from

$$\tilde{\sigma}_{eq} = \sqrt{(\tilde{\sigma}_1)^2 - \tilde{\sigma}_1 \tilde{\sigma}_2 + (\tilde{\sigma}_2)^2 + 3(\tilde{\tau})^2} \tag{23}$$

From the ratio of the yield stress $\sigma^Y$ to $\tilde{\sigma}_{eq}$, the strain vector at the yield point $\boldsymbol{\varepsilon}^Y$ is computed as

$$\boldsymbol{\varepsilon}^Y = \frac{\sigma^Y}{\tilde{\sigma}_{eq}} \boldsymbol{\varepsilon} \tag{24}$$

Therefore, the stress vector is obtained using $E$, $H$, $\mathbf{D}_1$, and $\boldsymbol{\varepsilon}^Y$ as

$$\begin{aligned}\boldsymbol{\sigma} &= \mathbf{D}_1 \boldsymbol{\varepsilon}^Y + \mathbf{D}_2 (\boldsymbol{\varepsilon} - \boldsymbol{\varepsilon}^Y) \\ &= \left(1 - \frac{H}{E}\right)\mathbf{D}_1 \boldsymbol{\varepsilon}^Y + \frac{H}{E}\mathbf{D}_1 \boldsymbol{\varepsilon}\end{aligned} \tag{25}$$

and the strain energy is computed from

$$S = \left(1 - \frac{H}{E}\right)\boldsymbol{\varepsilon}^T \mathbf{D}_1 \boldsymbol{\varepsilon}^Y + \frac{H}{2E}\boldsymbol{\varepsilon}^T \mathbf{D}_1 \boldsymbol{\varepsilon} \tag{26}$$

Although the stiffness after yielding depends on the stress ratio between $\sigma_{k1}$ and $\sigma_{k2}$, we assume the ideal state satisfying $\sigma_{k1} = \sigma_{k2}$ for which the relation between the equivalent stress and equivalent strain is obtained by experiment. Since we consider a monotonic loading process, the equilibrium shape of an air-pressured ETFE film can be obtained by minimizing $\Pi(\mathbf{X}) = S(\mathbf{X}) - W(\mathbf{X})$ with

$$S(\mathbf{X}) = \frac{1}{2}\sum_{k=1}^{m} A_k \left[ \boldsymbol{\varepsilon}_k^{YT} \boldsymbol{\sigma}_k^Y + \left(\boldsymbol{\varepsilon}_k - \boldsymbol{\varepsilon}_k^Y\right)^T \left(\boldsymbol{\sigma}_k + \boldsymbol{\sigma}_k^Y\right) \right] \tag{27}$$

# 5. Numerical examples

The proposed algorithm of approximate cutting pattern optimization is applied to a frame-supported PVC membrane (PVC-coated polyester fabric) and an air-pressured ETFE film. The optimization problems are solved using sequential quadratic programming implemented in SNOPT Ver. 7 [26], where the sensitivity coefficients are computed analytically by differentiating the governing equations. Triangular mesh is generated using the Delaunay triangulation library available in scipy.spatial package.

## 5.1 Model 1: Frame-supported PVC membrane

Consider an HP-type frame-supported membrane (Model 1) as shown in Fig. 5. Proportion of the model is $W_1 = 1.0W$, $W_2 = 1.3W$, and $H = 0.2W$. The material property is assumed to be orthotropic elastic as described in Sec. 2.1. Young's modulus in warp and weft directions are $E_x = 2.43 \times 10^2$ kN/m and $E_y = 2.27 \times 10^2$ kN/m, respectively. The shear modulus $G$ is 24.2 kN/m, and Poisson's ratios are $\nu_{xy} = 0.51$ and $\nu_{yx} = 0.55$. These material parameters are obtained from bi-axial and sheer tests according to the testing methods of Membrane Structures Association of Japan [27, 28]. The membrane is divided into two cutting sheets as shown in Fig. 6(a). The total numbers of nodes and elements are 160 and 240, respectively. The curved surface is projected to a plane that is parallel to the global $XY$-plane.

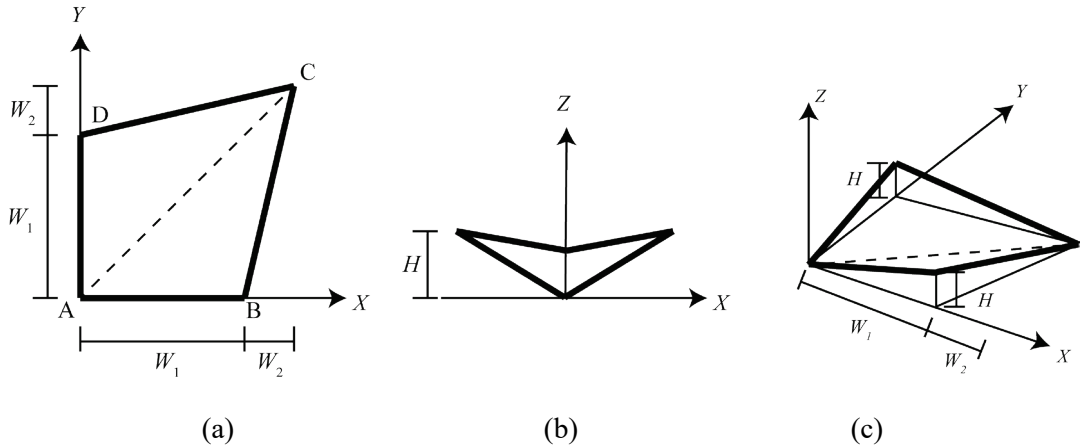

Figure 5: An HP-type frame-supported membrane structure (Model 1); (a) plan, (b) elevation, (c) diagonal view.

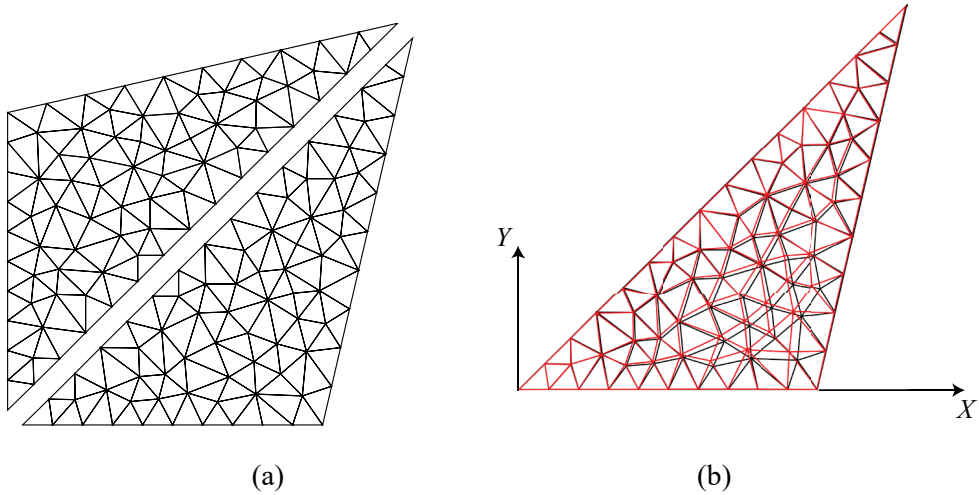

Figure 6: Cutting sheets of Model 1: (a) triangular mesh, (b) triangular mesh projected to *XY*-plane before optimization (black) and cutting sheet after optimization (red).

The target stress is 3.0 kN/m in both warp and weft directions. We first carried out the proposed process with different values of parameter $c$ for adjusting the convergence property in Eq. (20). As seen in Fig. 7, a small value of $c$ leads to gradual convergence, while a larger value results in rapid convergence after a large error in the first few steps. Based on this result, the parameter value $c = 0.5$ is chosen for Model 1.

The histories of average, maximum, minimum values, and standard deviation of stress are listed in Table 1, where *x* and *y* denote the directions of warp and weft, respectively. As seen from the table, the minimum value increases from a negative value to a positive value. The average value gradually converges to the target value. The minimum and maximum stresses are close enough to the target value at the 10th step. If we stop at the 20th step, the cutting pattern is as shown in Fig. 6(b). Note that the cutting pattern is close to the triangular plan of the half part of surface, which means that the area of cutting sheet is smaller than the surface area. The stress distribution at the 20th step is shown in Fig. 8. As seen from the figure, the stresses in warp and weft directions are almost uniform except in the area near corners. It is notable that almost uniform stress distribution can be achieved using only two cutting sheets.

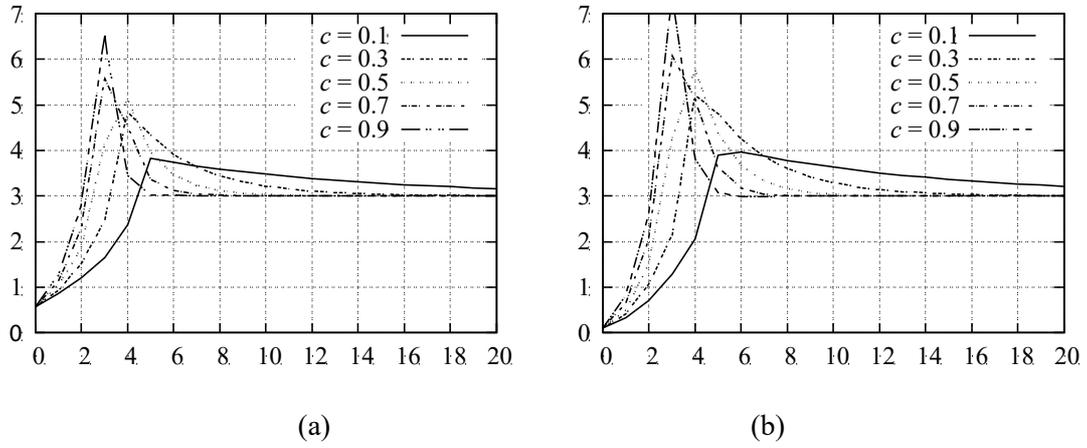

Figure 7: Variation of mean stresses for different values of parameter $c$ for adjusting the convergence property in Eq. (20); (a) warp direction, (b) weft direction.

Table 1: Histories of average, maximum, minimum values, and standard deviation of stress (kN/m) of Model 1.

|  | Step 0 | | Step 5 | | Step 10 | | Step 15 | | Step 20 | |
| --- | --- | --- | --- | --- | --- | --- | --- | --- | --- | --- |
| Direction | $x$ | $y$ | $x$ | $y$ | $x$ | $y$ | $x$ | $y$ | $x$ | $y$ |
| Average | 0.577 | 0.102 | 3.964 | 4.378 | 3.034 | 3.036 | 3.003 | 3.002 | 3.002 | 3.004 |
| Max. | 3.320 | 2.811 | 7.468 | 5.513 | 3.312 | 3.093 | 3.222 | 3.072 | 3.237 | 3.072 |
| Min. | −16.619 | −2.680 | 3.286 | 3.220 | 2.960 | 2.987 | 2.931 | 2.951 | 2.931 | 2.952 |
| Std. Dev. | 2.202 | 0.621 | 0.784 | 0.624 | 0.063 | 0.018 | 0.053 | 0.019 | 0.053 | 0.019 |

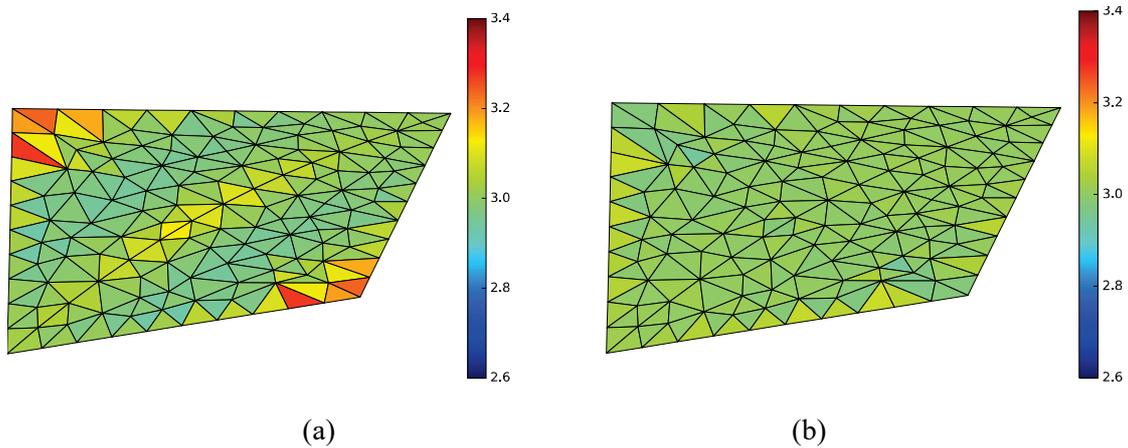

Figure 8: Stress distribution after 20 steps of optimization: (a) warp direction, (b) weft direction.

## 5.6 Model 2: Air-pressured ETFE film

Consider an air-pressured ETFE cushion, which has a square plan as shown in Fig. 9, where the ratio of $H$ to $W$ is 0.058. The material property is isotropic nonlinear elastic as described in Sec. 4. Thickness is 200 μm, Young's modulus is $1.60 \times 10^2$ kN/m, hardening coefficient is 10.4 kN/m, elastic shear modulus is 55.2 kN/m, shear modulus after yielding is 3.60 kN/m, Poisson's ratio is 0.45, and the stress and strain at yielding is 3.2 kN/m and 0.02, respectively. The specified air pressure is 1.0 kN/m$^2$, and the target stress is 4.0 kN/m. In this case, the radius of curvature is $2 \times 4.0/1.0 = 8.0$ m, if the surface is spherical. The curved surface is projected in the direction of the center of curvature of the surface as illustrated in Fig. 2(b).

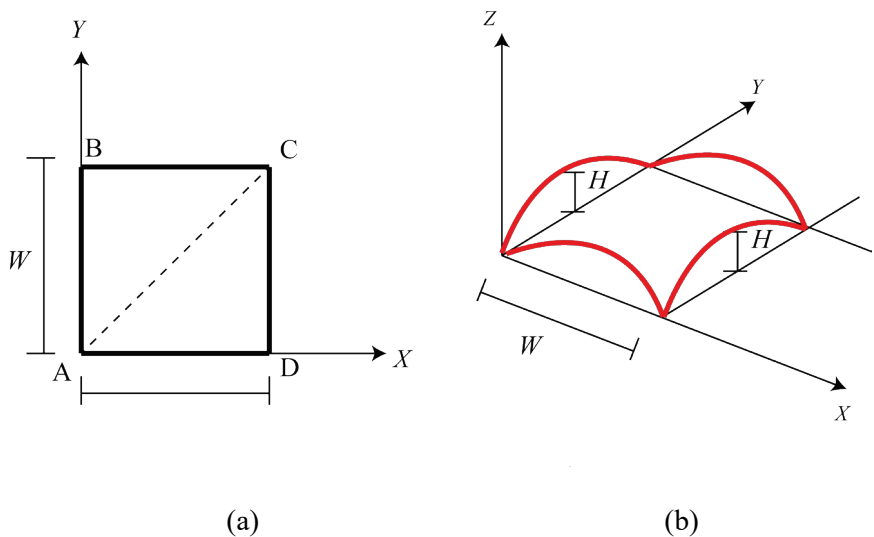

(a)  (b)

Figure 9: Air-supported ETFE sheet (Model 2); (a) plane, (b) diagonal view.

Table 2: Histories of average, maximum, minimum values, and standard deviation of stress (kN/m) of Model 2.

|  | Step 0 | | Step 2 | | Step 4 | | Step 7 | | Step 10 | |
| --- | --- | --- | --- | --- | --- | --- | --- | --- | --- | --- |
| Direction | $X$ | $Y$ | $X$ | $Y$ | $X$ | $Y$ | $X$ | $Y$ | $X$ | $Y$ |
| Average | 4.589 | 4.591 | 5.015 | 5.025 | 5.056 | 5.052 | 4.024 | 4.023 | 4.063 | 4.063 |
| Max. | 5.175 | 5.151 | 6.635 | 6.599 | 7.320 | 7.221 | 4.526 | 4.609 | 4.652 | 4.795 |
| Min. | -1.155 | -1.110 | 3.025 | 3.237 | 3.250 | 3.568 | 2.969 | 3.176 | 2.744 | 2.964 |
| Std. Dev. | 0.633 | 0.624 | 0.649 | 0.644 | 0.806 | 0.814 | 0.268 | 0.272 | 0.290 | 0.301 |

We carried out several trials for determining the value of $c$ in Eq. (20). Since the stress-strain relation is discontinuous for ETFE, a small value 0.05 has been chosen for $c$. The history of average, maximum, minimum values and standard deviation of stress is listed in Table 2, where $X$ and $Y$ denote the directions in $X$- and $Y$-directions on the global coordinates of the cutting sheet. As seen from the table, ETFE has a better accuracy than PVC, because the stiffness at the target stress of ETFE after yielding is smaller than that of PVC. The cutting pattern and stress distribution in $X$-direction after 10 steps are shown in Figs. 10(a) and (b), respectively. It is seen from these results that the cutting pattern is significantly different from the triangular shape, because the ratio of maximum height to the span $W$ has a moderately large value 0.2204. It is notable also for this model that almost uniform stress distribution can be achieved using only two cutting sheets.

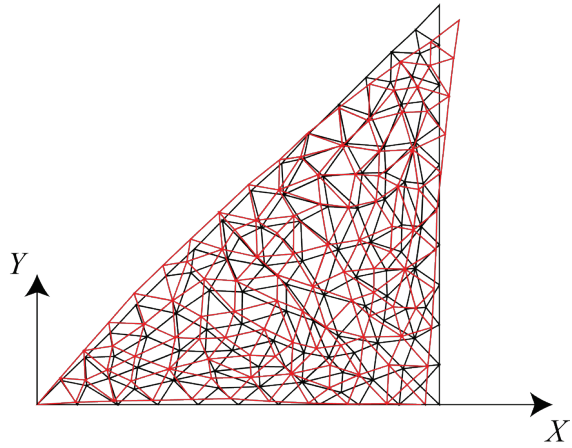

Figure 10: Cutting sheet of Model 2: triangular mesh projected to $XY$-plane before optimization (black) and cutting sheet after optimization (red).

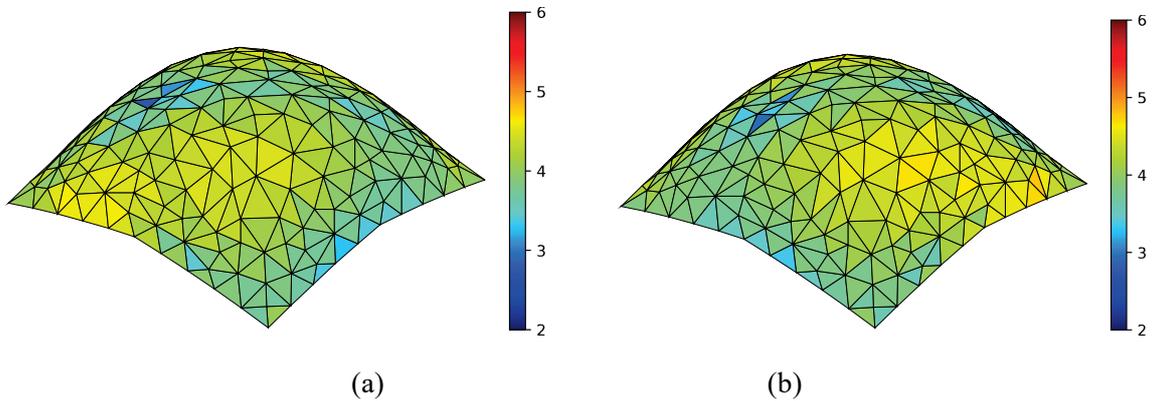

(a)  (b)

Figure 11: Stress distribution after 10 steps of optimization: (a) $X$-direction of initial shape; (b) $Y$-direction of initial shape.

*16*

## 6. Conclusions

An approximate method has been presented for cutting pattern optimization of membrane structures. The method is based on an inverse process of removing the specified target stress from the self-equilibrium state of curved surface. The surface and cutting sheets are divided into triangular elements with uniform plane stress, and the material properties are assumed to be isotropic elastic for the PVC sheet and biliniear nonlinear isotropic elastic for the ETFE film.

The approximate plane cutting pattern for the curved surface with specified target stress is obtained by removing the stresses of triangular elements and minimizing the norm of errors of the edge lengths from those of the triangular elements on a plane. By adjusting the stress parameter called reduction stress, an approximate optimal cutting pattern can be obtained after several iterations of cutting pattern generation and equilibrium shape analysis, which is formulated as an optimization problem of minimizing the strain energy under forced displacements at the boundary. Note that the reduction stress is regarded as a parameter for controlling the shape of cutting pattern; i.e., it does not represent the ideal target stress.

The equilibrium shape of a pneumatic membrane structure can also be obtained by minimizing the total potential energy including the work done by the air pressure. A formulation of pressure potential is obtained based on the shape variation in the normal direction of the surface. The material property of an ETFE sheet can be modeled as bilinear nonlinear elastic in the process of monotonically increasing the pressure to form the equilibrium shape. It has been shown in the numerical examples that almost uniform stress distribution is achieved for a frame-supported PVC sheet and an ETFE film subjected to air pressure using only two triangular cutting sheets.

## Acknowledgements

This study is partially supported by JSPS KAKENHI No. 16K14338.